\documentclass[aps,reprint,amsmath,amssymb]{revtex4}
\usepackage{graphicx}
\usepackage{dcolumn}
\usepackage{bm}
\usepackage{subfig}
\usepackage{amsmath,amssymb}
\usepackage{epstopdf}
\usepackage{textcomp}
\usepackage{hyperref}
\usepackage{listings}

\newcommand{\Fref}[1]{Fig.~\ref{#1}}
\newcommand{\Eqref}[1]{Eq.~(\ref{#1})}

\newcommand{\nn}{\nonumber}
\newcommand{\be}{\begin{equation}}
\newcommand{\ee}{\end{equation}}
\newcommand{\bear}{\begin{eqnarray}}
\newcommand{\eear}{\end{eqnarray}}

\newcommand{\ve}{\varepsilon}
\newcommand{\vem}{\varepsilon_\mathrm{mach}}
\newcommand{\ver}{\varepsilon_\mathrm{real}}
\newcommand{\vep}{\varepsilon_\mathrm{emp}}

\begin{document}

\title{Robust formula for $N$-point  Pad\'e approximant calculation
based on Wynn identity}
\author{Todor~M.~Mishonov}
\email[E-mail: ]{mishonov@bgphysics.eu}
\affiliation{Institute of Solid State Physics, Bulgarian Academy of Sciences,
72 Tzarigradsko Chaussee Blvd., BG-1784 Sofia, Bulgaria}
\affiliation{Faculty of Physics, St.~Clement of Ohrid University at Sofia,
5 James Bourchier Blvd., BG-1164 Sofia, Bulgaria}

\author{Albert~M.~Varonov}
\email[E-mail: ]{varonov@issp.bas.bg}
\affiliation{Institute of Solid State Physics, Bulgarian Academy of Sciences,
72 Tzarigradsko Chaussee Blvd., BG-1784 Sofia, Bulgaria}
\affiliation{Faculty of Physics, St.~Clement of Ohrid University at Sofia,
5 James Bourchier Blvd., BG-1164 Sofia, Bulgaria}

\date{21 April 2020}

\begin{abstract}
The performed numerical analysis reveals that Wynn's identity for the compass
$1/(N-C)+1/(S-C)=1/(W-C)+1/(E-C)=1/\eta$
(here C stands for center, the other letters correspond to the four directions of the compass)
gives the long sought criterion, the minimal $|\eta|$, for the choice of the optimal Pad\'e approximant.
The work of this method is illustrated by calculation of multipoint Pad\'e approximation by
a new formula for calculation of this best rational approximation.
The work of the criterion for the calculation of optimal Pad\'e approximant is illustrated by the frequently seen in the theoretical physics problems of calculation of series summation and multipoint Pad\'e approximation used as a predictor for solution of differential equations motivated by the magneto-hydrodynamic problem of heating of solar corona by Alv\'en waves.
In such a way, an efficient and valuable control mechanism for $N$-point Pad\'e approximant calculation is proposed.
We believe that the suggested method and criterion can be useful for many
applied problems in numerous areas not only in physics but in any scientific application where differential equations are solved.
The solution of the Cauchy-Jacobi problem is illustrated by a Fortran program.
The algorithm is generalized for the case of the first $K$ derivatives at $N$ nodal points.
\end{abstract}

\maketitle

\section{Introduction and motivation}

As the analytical structures take central place in the theoretical and mathematical physics,
the Pad\'e approximants of functions become one of the most important problems of applied mathematics.
For a general introduction in the problem of Pad\'e approximants, see the well-known monographs Refs.~\cite{Pade,Brezinski:91,Brezinski:96}.
Since the problem was first systematically studied by Frobenius~\cite{Brezinski:H,Frobenius:81} and the Pad\'e table introduced by Pad\'e~\cite{Brezinski:H,Pade:92} in the 19$^\mathrm{th}$ century,
the contemporary literature is enormous.
We will mention only two popular problems, which can be of broad interest for physicists:
1) series summation and 2) extrapolation of functions related to predictor-corrector methods for solution of ordinary differential equations.
Searching for numerical recipes in this field, a beginner can read that the downside of the Pad\'e approximation is that it is uncontrolled and the dangers of extrapolation cannot be overemphasized~\cite{NumRep}.
And after that our beginner will most probably start searching for a commercial software, where the appropriate analytical formulae are programmed.
By the very famous law: whatever can go wrong, will go wrong, in very interesting for the physics problems, the commercial software stops working and the research would continue only in co-authorship with the source code authors.

In 1947 Hannes Alvf\'en~\cite{Alfven:47} concluded that the solar corona heating is due to the absorption of magneto-hydrodynamic (now called Alfv\'en) waves.
In order to describe the almost 100-time increase of the plasma temperature in a very narrow region of the solar atmosphere.
The programming of the magneto-hydrodynamic equations is relatively easy, but after 30-time temperature increase all available for us commercial computation software stops working.
Facing such a problem, we had to develop an own method for a frequently encountered numerical problem (four first order linear differential equations coupled to four first order non-linear differential equations)~\cite{MHD,MHD:19}.
Our research goal was to develop this method based on calculation of Pad\'e approximants with a mechanism to select the most accurate approximant.
This mechanism is precisely the main development of the current research and allows the usage of the described method for solution of differential equations.
That is why we suppose that many colleagues are in our situation and we share know-how in calculation of Pad\'e approximants.

Before starting let us recall the definition of multipoint Pad\'e approximation\cite{Pade}. 
A rational function $f(z)$ which fits given function values $y_i$ 
at various points $x_i$, $i=1,\dots,\, N$, i.e. $f(x_i)=y_i$ is called multipoint Pad\'e approximant.
The associate problem of interpolation by rational functions is called Cauchy-Jacobi problem.
Multipoint Pad\'e approximants are also called:
rational interpolants,
N-point Pad\'e approximants,
Newton-Pad\'e approximants, etc,
depending on the context. 
Addressing to experts 
in the present paper we present a
general robust solution of the Cauchy-Jacobi problem
which can attract interest and become a convenient tool in the applied numerical analysis.

\section{Announcement of the results}

For a practical implementation of an analytical result, we have to use a computer with finite accuracy of digital representation of real numbers.
The Pad\'e approximants give fast convergence and actually the best rational approximation 
but they are sensitive with respect to the noise of the discrete representation of real numbers.
Roughly speaking every value representing a real number suffers from error of discretization 
as it is multiplied by a random factor
$1+\vem (\mathrm{rnd}-1/2)$,
where $\vem$ is the machine epsilon (the biggest positive value for which $1+\vem=1$),
and ``$\mathrm{rnd}$'' is a programming operator generating a random number
homogeneously distributed  between 0 and 1.

As formulas for Pad\'e approximants picturally speaking are very sensitive
to the random noise of truncation,
it is necessary to study the influence of the machine error $\vem$ on the final result.
In such a way beyond the analytical problem 
we have also a statistical one
which already belongs to another branch of mathematics.
Probably for some special cases it is possible to calculate the probability distribution functions (PDF) of the errors of the final result of calculation of Pad\'e approximants, 
but it deserves to start with descriptive statistics of some well-known examples, which in our opinion illustrate how to calculate and choose the optimal Pad\'e approximant in physical applications.
Here the determination of the optimal Pad\'e approximant is crucial, as our beginner has already learned from the numerical recipes book.
To experts we propose a robust method which
at least empirically works well.

Let $r_{l,m}(z)$ are different Pad\'e approximants of the complex function $f(z)$
\be
f(z) \approx r_{l,m}(z) = \frac{\sum\limits_{i=0}^l a_i z^i}{\sum\limits_{k=0}^m b_k z^k},
\qquad b_0=1,
\ee
where $l$ is the maximal power of the numerator and $m$ is the maximal power in the polynomial
in the denominator.
For these approximants in 1966 Wynn discovered
the remarkable relation~\cite[Eqs.~(15) and (16)]{Wynn:66}
\be
\frac{1}{\eta_{lm}} \equiv \frac1{r_{l,m+1}-r_{l,m}}+\frac1{r_{l,m-1}-r_{l,m}}
=\frac1{r_{l+1,m}-r_{l,m}}+\frac1{r_{l-1,m}-r_{l,m}},
\label{Eq:Wynn}
\ee
later on baptized by Gragg~\cite[Theorem~5.5]{Gragg:72} as \textit{missing identity of Frobenius}.
The criterion for choosing of optimal Pad\'e approximant is based on this 
Wynn relation which can be rewritten also as
\begin{align}&
\frac{1}{\eta} \equiv \frac1{S-C}+\frac1{N-C}
=\frac1{E-C}+\frac1{W-C},\\&
C\equiv r_{l,m}, \quad
N\equiv r_{l,m-1}, \quad 
E\equiv r_{l+1,m},\quad
W\equiv r_{l-1,m},\quad
S\equiv r_{l,m+1}.
\label{Eq:Wynn_compass}
\end{align}
In the analysis of algorithms is is convenient to represent different
approximants $r_{l,m}$ in a Pad\'e table.
If $C$ is in the Center, the notations of the other elements $N,\,E,\,W,\,S$
follow the direction of the compass: North, East, Wests, and South, also found in the original paper of Wynn~\cite{Wynn:66} (i.e. we use Central NEWS algorithm).
Probably this representation was an idea of Gragg, who has also described it~\cite{Gragg:72}, because of the footnote in Wynn's paper on page 266 (if that's the case, Gragg was a referee of Wynn's paper).
For convergent Pad\'e approximants
\be
f(z)=\lim_{l,m\rightarrow\infty}r_{l,m}(z)
\ee
in the Pad\'e table $C,\,N,\,E,\,W,\,S\rightarrow f(z)$, $1/\eta\rightarrow \infty$,
and
\be
\lim_{l,m}\eta_{lm}\rightarrow 0.
\ee
In numerical implementations this limit leads to minimal-$|\eta|$ criterion at CNEWS algorithm.

Let us repeat in other words. 
Often the difference between sequential Pad\'e approximants gives a reasonable evaluation of the 
error, see Ref.~\cite{Pade}.
In case of convergence, we have vanishing differences between the values of different cells of the Pad\'e table $(r_{l,m+1}-r_{l,m}) \rightarrow 0$, for $l,m \rightarrow \infty$  and $l/m=\mathrm{const}$.
In this case also $\eta_{lm} \rightarrow 0$ and the minimal value of $\eta_{lm}$ gives a reasonable 
criterion for the minimal error and the choice of the optimal Pad\'e approximant. 
This theoretical hint we prove on many examples of calculation of Pad\'e approximation
and arrive at the conclusion that
the long sought criterion for the choice of the optimal Pad\'e approximant is simply a search for the minimal value $|\eta_{l,m}|$, i.e.
\be
\eta_\mathrm{min} \equiv |\eta_{L,M}|
=\mathrm{min}_{l,m} |\eta_{l,m}|
\ee
and $f(z)\approx r_{L,M}$.
It is technologically and aesthetically attracting that the Wynn identity gives simultaneously
1) explicit method for the calculation of Pad\'e approximants in $\ve$-table and 
2) method for the evaluation of the error; $\lim \eta_\mathrm{min}=0$ is actually a criterion for convergence of  Pad\'e approximation if we consider real numbers.
Our criterion is an alternative of the singular value decomposition (SVD) method described in great detail
in Refs.~\cite{Gonnet:13,Beckermann:13}, see also Refs.~\cite{Karlsson:76,Gonchar:78,Gonchar:89,Levis:18,Yattselev:18}. 
Our minimal $\eta$ criterion and the tolerance level of SVD 
are perhaps different implementations of one and the same idea.
For both methods the rounding error is something 
external for the theory of Pad\'e approximation. 
There are no theoretical justifications which method is better 
and the numerical experiment and the descriptive statistics are just the first steps 
in the practical realizations on the robust Pad\'e approximants.

The purpose of the present work is to demonstrate how this criterion works
and can be used in applications.
In the next section we illustrate in detail this criterion on the problem 
of of calculation of limes of a numerical series.

\section{The new algorithm for determination of the optimal Pad\'e approximant}

Let us have a numerical sequence $\{ S_0, S_1, S_2, \dots , S_N \}$ and we need to calculate the limit $S=\lim\limits_{l \rightarrow \infty} S_l$.
The well-known method is to initialize
$r_{l,0}=S_l$, for $l=0, \dots, N$, and also $r_{l,-1}=\infty$ for $l=0, \dots, N+1$.
In other words we treat every term of the initial series as
polynomial approximation of some function
\be
r_{l,0}=\sum\limits_{i=0}^l a_i z^i= S_l.
\ee
The method can work even if this sequence is convergent,
and those are the most interesting cases for applications.
Then we calculate in the ``south direction'' the corresponding Pad\'e approximants
\be
S=C+\dfrac{1}{\dfrac{1}{E-C}+\dfrac{1}{W-C}
-\dfrac{1}{N-C}}
\label{sought}
\ee
or
\be
r_{l,m+1}=r_{l,m}+\dfrac{1}{\dfrac{1}{r_{l+1,m}-r_{l,m}}+\dfrac{1}{r_{l-1,m}-r_{l,m}}
-\dfrac{1}{r_{l,m-1}-r_{l,m}}}
\label{sought}
\ee
and simultaneously calculate the empirical error
\be
\eta_{lm} = \dfrac{1}{\dfrac{1}{r_{l+1,m} - r_{l,m}}+\dfrac{1}{r_{l-1,m}-r_{l,m}}}
\equiv\eta= \dfrac{1}{\dfrac{1}{E - C}+\dfrac{1}{W-C}}
. \label{eta}
\ee
The minimal absolute value in the $\eta$-table $\eta_\mathrm{min}$ is our criterion for the determination of the optimal Pad\'e approximant.
According to the best we know, this 
minimal-$|\eta|$ criterion has never been implemented in the numerical recipes so far (authors will very much appreciate any information about this).

For the programming task, we have to calculate all the values, for which division is possible.
In the next section we thoroughly describe 2 technical examples.

\section{Two easy pieces}

\subsection{Calculation of divergent series $\ln(1+x)$}

\begin{figure}[h]
\includegraphics[scale=0.4]{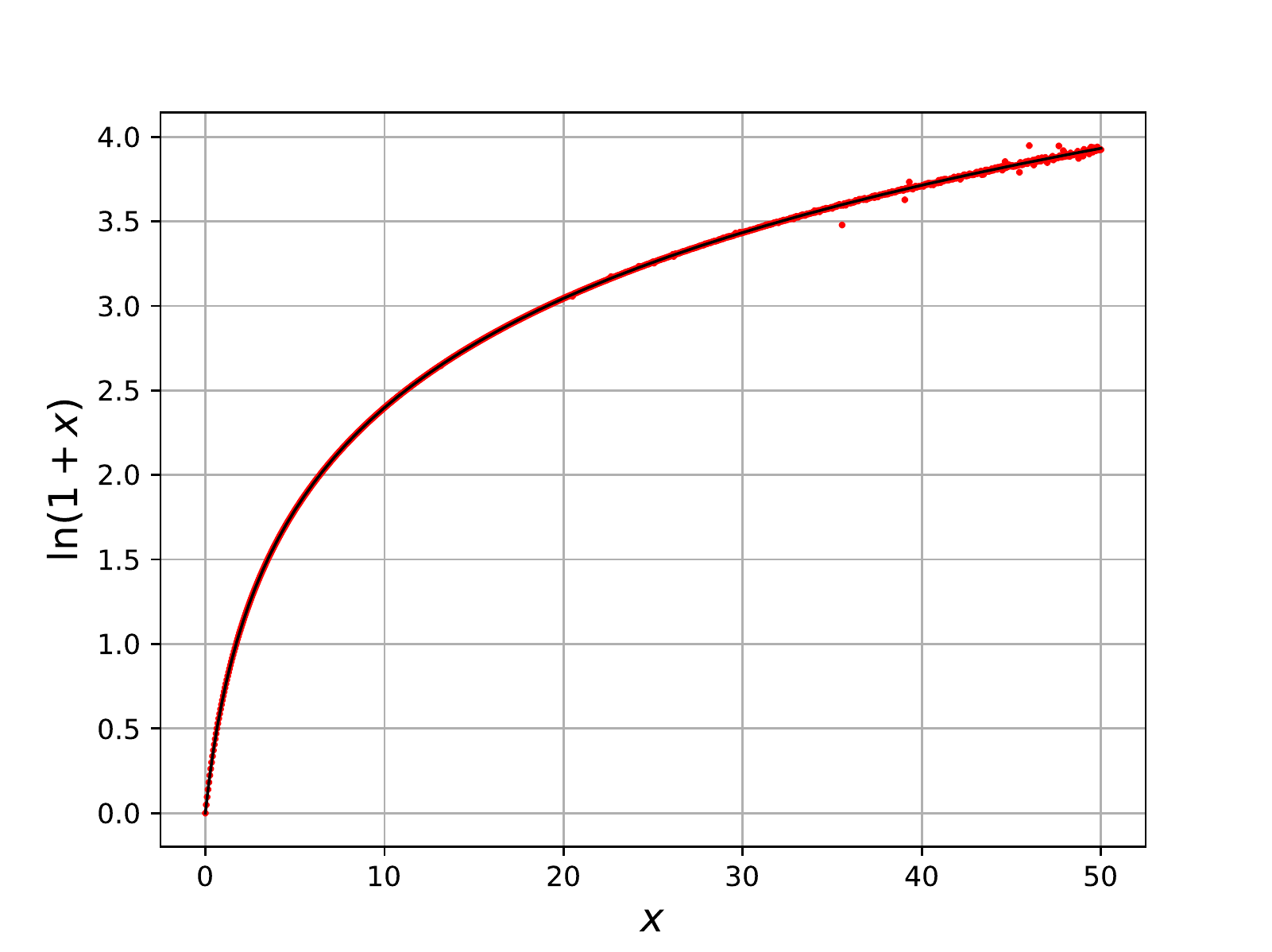}
\caption{The logarithmic function (line) and the series summation 
with $\eta_\mathrm{min}$ criterion for the optimal choice (dots). 
The calculated values are \textit{evaporated} from the analytical curve, but $\eta_\mathrm{min}$ criterion gives reliable warning depicted in \Fref{fig:log101ee}.}
\label{fig:log101}
\end{figure}

\begin{figure}[h]
\includegraphics[scale=0.4]{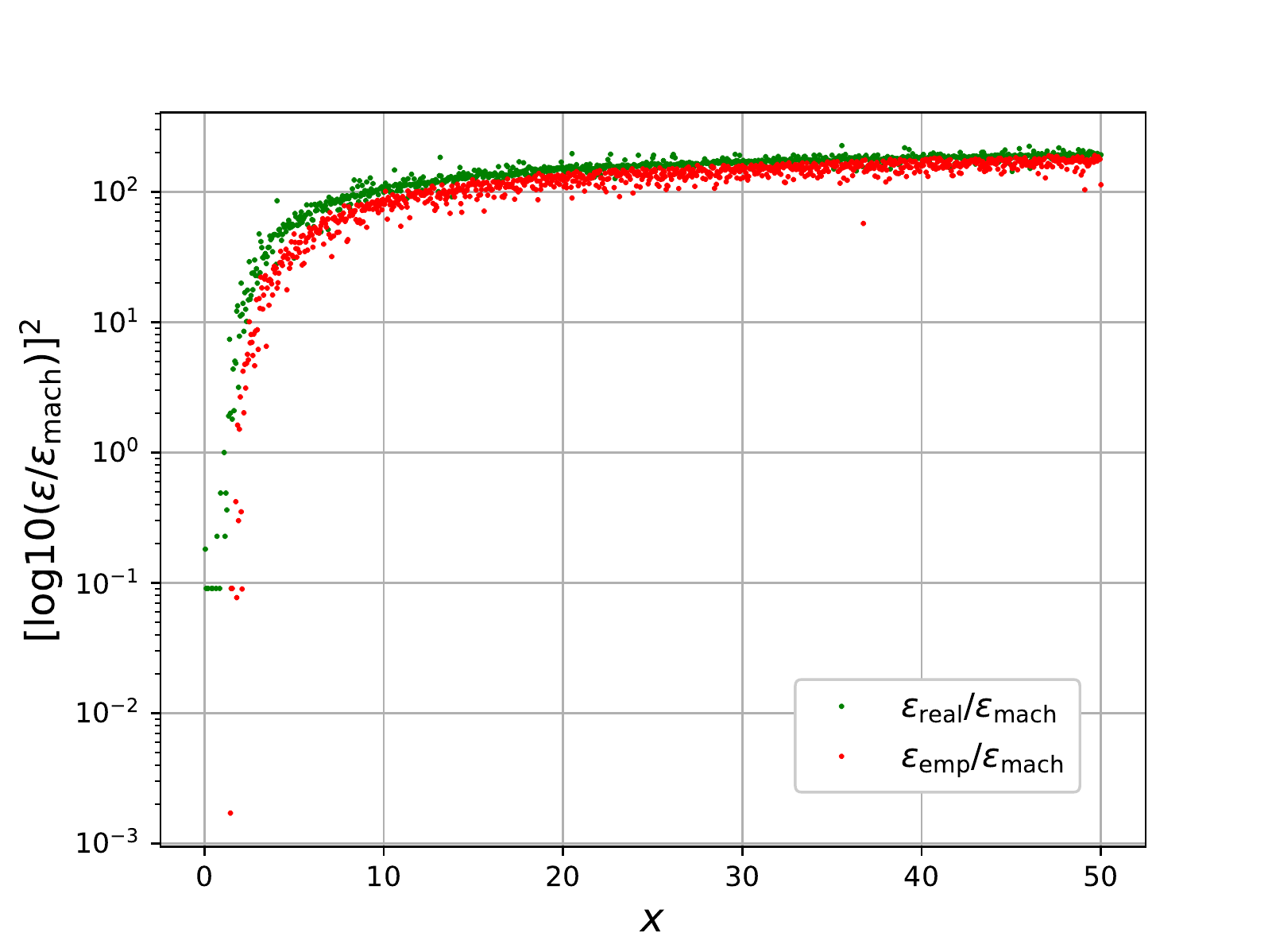}
\caption{Square of decimal logarithms of empirical $\vep \equiv \eta_\mathrm{min}$ and real $\ver$ error versus the argument of the function for the calculation of the $\ln(1+x)$ series in \Fref{fig:log101}.
Close to the convergence radius the errors are small and almost linearly increase,
then the errors reach saturation when the numerical resources of the fixed accuracy are exhausted.}
\label{fig:log101ee}
\end{figure}

\begin{figure}[h]
\includegraphics[scale=0.4]{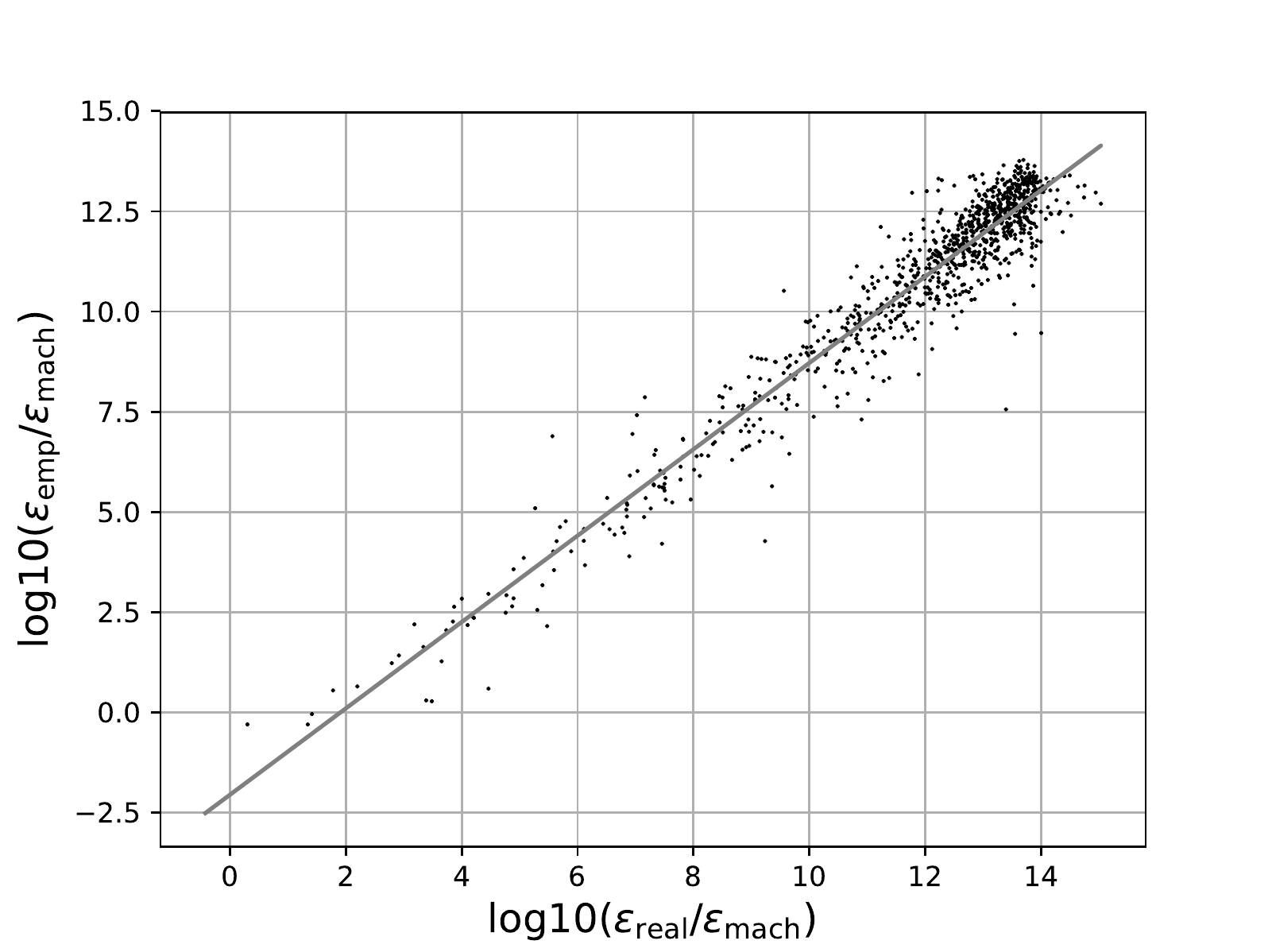}
\caption{The logarithm of the empirical $\eta_\mathrm{min} \equiv \vep$ versus logarithm of the real $\ver$ error for the calculation of the $\ln(1+x)$ series in \Fref{fig:log101}.
The high correlation coefficient 0.961 of the linear regression reveals that the long sought criterion for empirical evaluation of the accuracy of the Pad\'e approximants calculated by $\ve$-algorithm has already been found.}
\label{fig:log101e}
\end{figure}

Perhaps the most famous example summation 
of divergent series by Pad\'e approximants is the calculation of the divergent Taylor series $\ln(1+x)$ 
beyond the radius of convergence $|x|<1$
\be
\ln(1+x)=\lim_{n\rightarrow\infty} S_n, \qquad 
S_n(x) \equiv \sum_{k=1}^n (-1)^{(k+1)} \frac{x^k}{k}.
\label{eq:log1x}
\ee
The application of Eq.~(\ref{sought}) and Eq.~(\ref{eta}) gives the optimal Pad\'e approximant for every positive $x$ shown in \Fref{fig:log101}.
Even for $x \approx 20$ for the calculation of the series $\ln(1+x)$ a pixel accuracy is present.
In many articles results of numerical calculations are illustrated by the figures
presented by computer graphics with million pixels.
For such a figures pixel accuracy means
maximal error in which the approximation is undistinguished from the exact 
result when graphically presented on screen. 
For many applied studies this is an acceptable beginning.
For larger values of $x$ the dots representing the calculation \textit{evaporate} from the line representing the exact value.
The most important detail is how and when the method stops working and when the resources of the numerical accuracy are exhausted.
The real $\ver$ and empirical $\vep \equiv \eta_\mathrm{min}$  of the calculation are shown in \Fref{fig:log101ee}.
The empirical error $\vep \equiv \eta_\mathrm{min}$ shows saturation and is a reliable indicator for the calculation accuracy.
Statistically analysis of the correlation between both errors from \Fref{fig:log101ee} shows very strong dependence between them in \Fref{fig:log101e}.
Now we can safely tell how accurate our Pad\'e approximation is and even how far out in $x$ it can be extended.
This dependence  reveals that the long sought criterion for empirical evaluation of the accuracy of the Pad\'e approximants calculated by $\ve$-algorithm has already been found.

\subsection{Extrapolation of a sine arch from a preceding one}

\begin{figure}[h]
\includegraphics[scale=0.4]{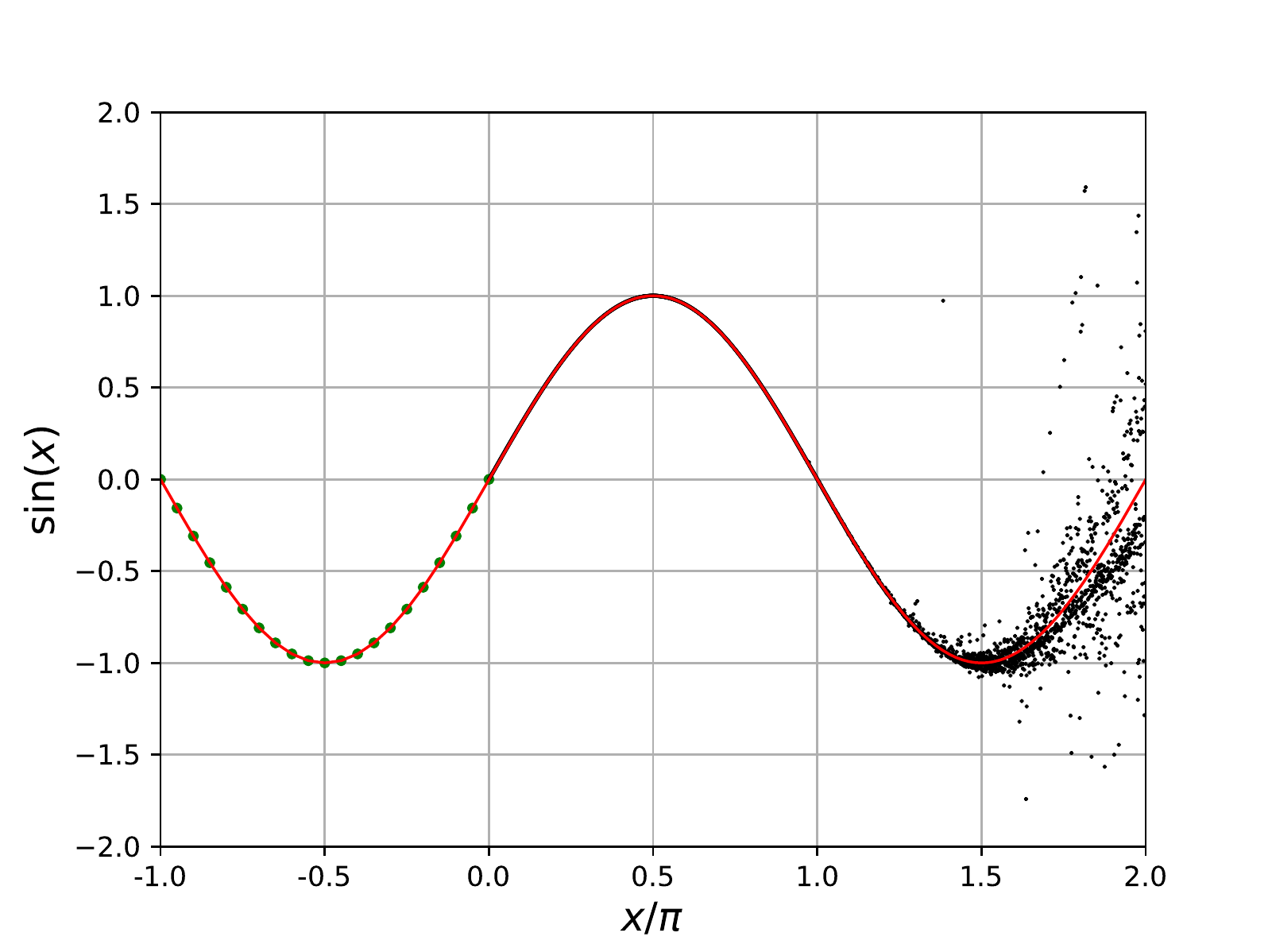}
\caption{Extrapolation of the function $\sin(x)$ (small dots) from
21 interpolation points (larger dots) in the interval $[-\pi,0]$ compared with the real function (line).
In the interval $(0,\pi)$ the function is reliably extrapolated and the limit of the numerical implementation of the Aitken-Wynn extrapolation is clearly shown -- a \textit{gas} of extrapolated points \textit{evaporated} from the analytical function.
}
\label{fig:s21}
\end{figure}

\begin{figure}[h]
\includegraphics[scale=0.4]{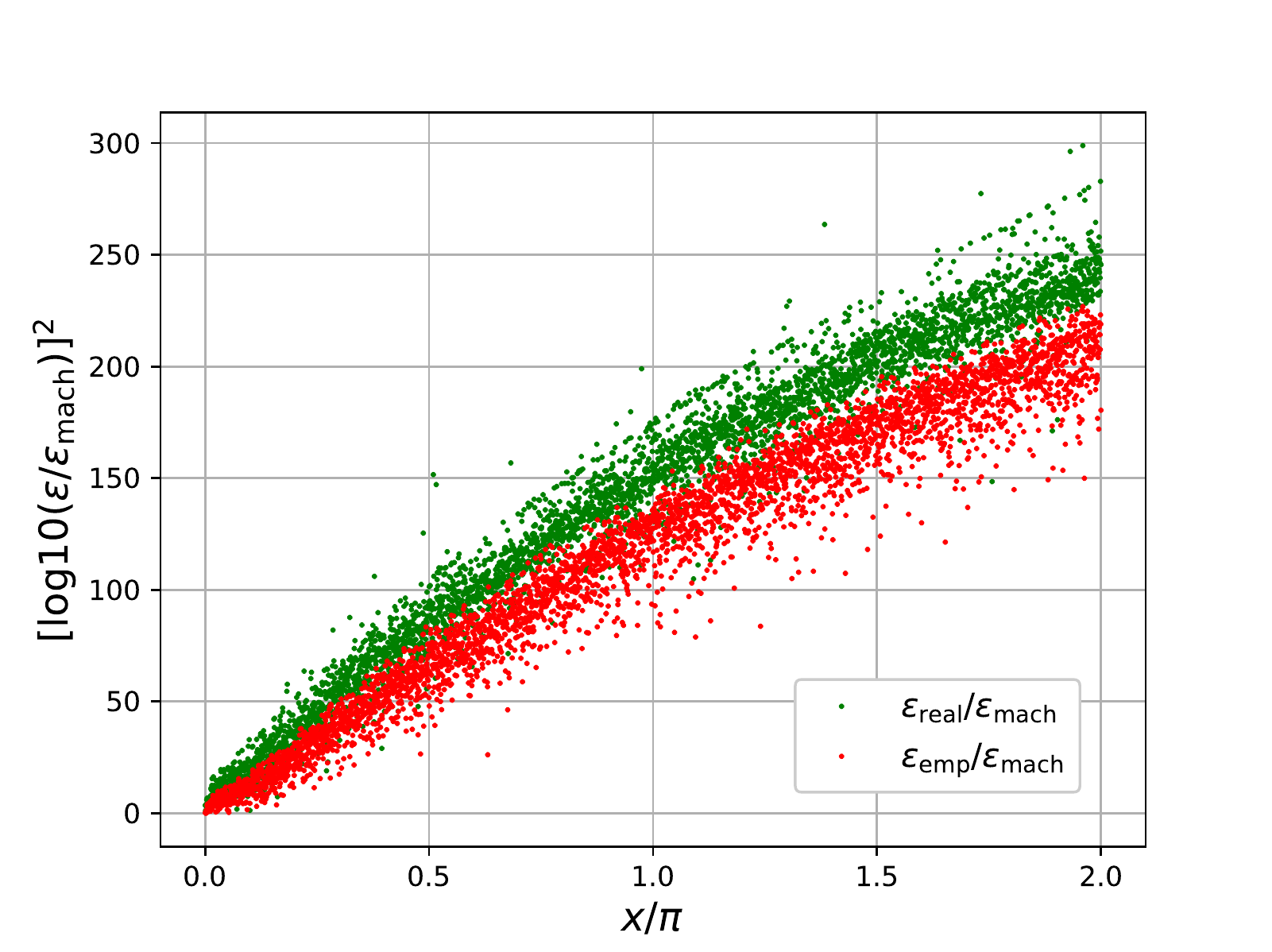}
\caption{Squared logarithm of the error estimates $\varepsilon$ of the $\sin(x)$ extrapolation shown in \Fref{fig:s21}.
The error $\ver$ is the real error of the extrapolation and the error $\vep \equiv \eta_\mathrm{min}$.
The similar behaviour of both errors shows that our criterion is a reliable method for error estimation, 
which is evident in \Fref{fig:s21e}.
The important problem in front of the applied mathematics is to research real extrapolation error beyond the extrapolation interval.
}
\label{fig:s21ee}
\end{figure}

\begin{figure}[h]
\includegraphics[scale=0.4]{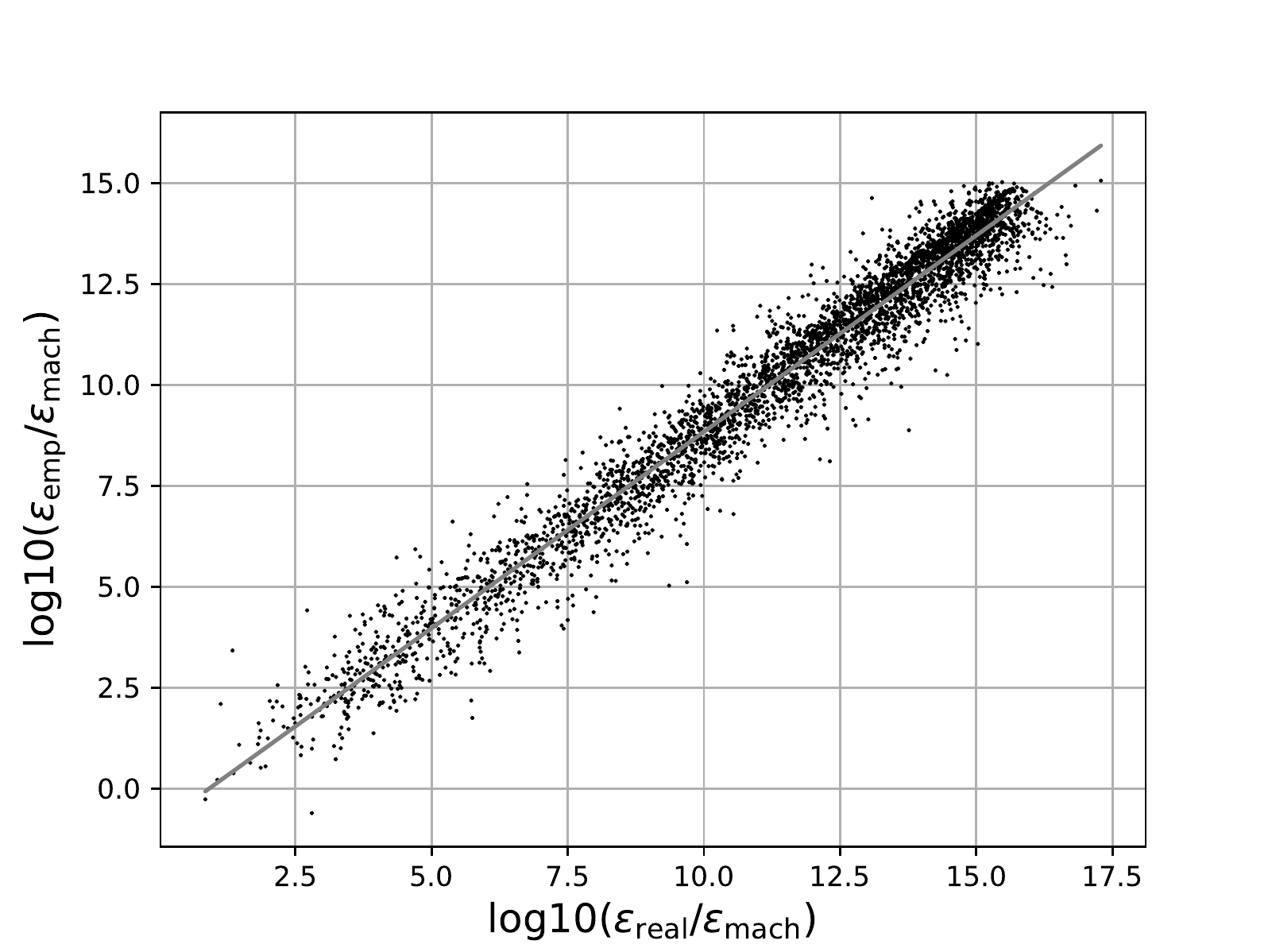}
\caption{Logarithmic dependence of $\eta_\mathrm{min} \equiv \vep$ as a function of the real error $\ver$ for the extrapolation shown in \Fref{fig:s21}.
The slope of the line of the linear regression is 0.973, the intercept is -0.887 with a correlation coefficient 0.979.
This strong correlation sets in the agenda the problem of the statistical properties of the Pad\'e approximants and the investigation of the corresponding probability distribution functions (PDF)}
\label{fig:s21e}
\end{figure}

The second technical example we present is the problem of extrapolation of function which we illustrate 
in the case of the $\sin(x)$ function.
We take $N$ equidistant interpolation points from one arch of the $\sin(x)$ function
and extrapolate the next arch and even beyond. 
We use the well-known Aitken interpolation method~\cite{Aitken:32,Bronshtein:55,Abramowitz:72,Korn2} to order the interpolation points within the arch and the point to be extrapolated.                                                                                                      
In more details:
we calculate $y_i=\sin(x_i)$ for all $N$ points $i=1,\dots,N.$                                                                                                    
The digital noise is minimized if the points are ordered
with the proximity with the point of interpolation, i.e.
\be
|x-x_1|\le |x-x_2|\le |x-x_2|\le\dots \le |x-x_N|.
\ee
Then using sequentially $l$ interpolation points in the Aitken scheme
we calculate the corresponding polynomial approximation $S_l$ for $l=1,\dots,N.$
Next we apply the CNEWS method and use the minimal-$|\eta|$ for the best Pad\'e
approximant.                                                                                                                                                                                                                                            
In this manner we have a numerical sequence that we give to the $\ve$-algorithm to calculate its limit.
The described calculation for the $\sin(x)$ function is shown in \Fref{fig:s21}.
The preceding arch is in the interval $(-\pi,0)$ and contains 21 interpolation points.
Using these points, we extrapolate one arch with 2000 points in the interval $(0,\pi)$ and continue the extrapolation in an attempt to obtain a second arch with the same number of points in the next interval $(\pi,2\pi)$.
Repeating, every extrapolated point is calculated independently using the Aitken method for interpolation followed by calculation of Pad\'e approximants by the Wynn identity and the optimal Pad\'e approximant is chosen by the advocated in the present paper $\eta_\mathrm{min}$ criterion.
The deviation of the points from the real function represented with the line in \Fref{fig:s21} shows the limit of applicability of the Aitken-Wynn extrapolation algorithm.
Detailed error estimates of the extrapolation are shown in \Fref{fig:s21ee}.
The $\eta_\mathrm{min}$ criterion shown in \Fref{fig:s21ee} gives the reliable order estimation of the error.
The similar behaviour of both errors shows that our criterion is a reliable method for error estimation.
Furthermore, a statistical analysis of the correlation between the real $\ver$ and the empirical error $\vep$ shows very strong dependence in \Fref{fig:s21e}.

The literature on the problem of extrapolation is enormous, however we have not found a discussion of the problem of the mathematical statistics.
If we know the values of the function with some accuracy, how to choose the best Pad\'e approximant among the many, some of which have serious noise of rounding and others are even \textit{uncontrolled}~\cite{NumRep}.
There is a huge literature about the formulae, but a convenient criterion for the choice of optimal Pad\'e approximant is not discussed.
We give our illustrative example not as an exercise of programming, but as an illustration that we have a good working criterion for many practical cases.
As a rule, articles on Pad\'e approximation are illustrated by examples describing how some method is good working, but the instructive analysis reveals where the method stops working and what is most important whether some criterion warns that calculation resource has been exhausted.
From a practical point of view, exhaustion of resources for extrapolation is demonstrated by \textit{evaporation} of the extrapolated points from the analytical curve, and one of the achievements of the present work is that order of the magnitude of this deviation is reliably indicated by our criterion based on the Wynn identity.  

\section{Discussion and conclusions}

In spite that the literature for Pad\'e approximants is enormous we were unable to find alternative formulas, methods, algorithms and programs to compare the work of our method.
One can put in the agenda for the development of the numerical analysis a simple test:
which method for extrapolation of functions from tabulated values in $N$ points gives the best extrapolation of a function far from the interpolation interval.
This juxtaposition will give the final verdict which method is appropriate to be implemented in commercial  software like Mathematica and Maple. 

The new result of the implementation of calculation of Pad\'e approximants and their application is the modulus minimization $\eta_\mathrm{min}$ of Wynn's $\eta_{lm}$ as reliable empirical criterion of the error.

The preformed analysis of several simple examples has revealed that for practical implementation of 
Pad\'e approximants we can use the empirical $\eta_\mathrm{min}$ error extracted from the Wynn's identity. 
In the agenda the statistical problem of calculation of PDF of the Pad\'e approximants has already been set. 
The comparison of descriptive statistics data for the PDF of errors of calculation of Pad\'e approximants by different criteria will give what the answer what general recommendation as a numerical recipe have to be given to users not willing to understand how.

In short, the practical implementation of Pad\'e approximants can reach one order of magnitude more applications in theoretical physics and applied mathematics.
More than half a century after its discovery, the $\ve$-algorithm has not been included for calculation of divergent series with convergent Pad\'e approximants and for extrapolation of functions in commercial software.
Now the time for this inclusion has come, the herein implemented control mechanism has rendered this mission possible.

Last but not least, the suggested criterion \Eqref{eta} is applicable in solution of differential equations, numerical analytical continuation, perturbation, series summation and other analogous problems of theoretical physics.

Let us summarize our results: 
1) we have suggested a method for empirical choice of optimal Pad\'e approximant which woks for all examples we have met in the literature and seems to have universal applicability,
2) we derived a new formula for the $N$-point Pad\'e approximant, for a function tabulated in $N$-points.
This formula based on Eitken interpolation formula and Wynn identity is maximally robust with respect of truncation errors in the numerical calculations and can be used even for extrapolation.
We have developed this method in order to obtain one predictor method for predictor-corrector methods for solving of one irrelevant for the present paper magnetohydrodynamic problem,
the programs where our formulae are implemented are given as appendices of the arXiv version of the present paper~\cite{arxiv} and early and more detailed version has already been published in conference proceedings~\cite{AIP_Pade}.
The figures in the paper are calculated by those programs.
Reproducing the figures is the criterion that our formula can be used directly in numerical calculations.
Up to now a non-specialist knows that extrapolation and summation of divergent series is a forbidden procedure used only by elite mathematicians. 
We have made a popularization addressed to user of applied mathematics working for the industry.
3) The described method can be used as a predictor method for solving
ordinary differential equations $\mathrm{d}y/\mathrm{d}x=\mathcal{F}(x,y)$.
The described N-point Pad\'e approximation can be used as a predictor
using former calculated points and their first derivatives as interpolation points
to calculate the predicted value of the function $y_\mathrm{pred}$
at the new point $x_\mathrm{new}$ of the argument.
Then we calculate the derivative at the new point 
$\mathcal{F}(x_\mathrm{new},y_\mathrm{pred})$.
We consider this possibility as a new method for solution of ordinary differential equations.
The method can be additionally improved if in the ($x$-$y$) plane we perform a rotation,
and the new abscissa of the interpolation to be taken along the local tangent of the curve.
We started with the problem of calculation of the temperature profile of the solar corona
and velocity of solar wind\cite{MHD:19} but this method can be applied for many other problems for which standard software is not (good) working.
4) If we know the first $K$ derivatives of the function using of Taylor expansion 
we can calculate
\be
\label{Taylor}
\tilde{y}_i=y_i+(x-x_i)\,y_i^{\prime}+\frac12\,(x-x_i)^2y_i^{\prime\prime}+\dots
+\frac1{K!} (x-x_i)^K\,y_i^{(K)}.
\ee
Then we can use the Taylor approximants $\tilde{y}_i$ as interpolation points $x_i$ in the already described method.
At the basic points the interpolated rational function will have exact coincidence not only for the function 
but also for its first $K$ derivatives.
The main idea of the algorithm is the same:
using a node point $x_i$ we calculate the interpolated value 
$\tilde{y}(x)$ at the argument of interpolation, and then
use this calculated value $\tilde{y}(x;x_i)$ at this point as a new value
in the node point $y_i=f(x_i)\leftarrow \tilde{y}(x;x_i)$.
If we need a formal proof of convergence it is necessary in the right
side of Eq.~(\ref{Taylor}) to put insert the evaluation of error
$h^{(K+1)}M_{(K+1)}/(K+1)!$, where $h$ is the maximal distance 
between interpolating an node argument and $M_{(K+1)}$
is the supremum of the modulus of the $(K+1)$ derivative.
This estimation of the error reveal why the analytical approximation is 
not applicable to function $|x|$ and the trajectory of a Brownial particle.

Finally we have to remark. 
There is no general criterion for the convergence of the Pad\'e approximation,
and even to choose the optimal approximant when the method is convergent,
that is why there are no theorems which could be cited.
As all known existing methods are given without any empirical prescription 
for the choice of the optimal approximant the considered minimal-$|\eta|$ criterion
of the CNEWS method can not be compared with other methods.
The area of the applicability of the suggested method is however well defined.
For many applied problems it is known that solution exists, while
for the analytical problems the solution is analytical.
In many problems in the physics, for example,
there are alternative method to check whether the searched solution is numerically obtained
in spite that there will be no alternative methods.
Contra-examples of mathematical series
for which the Pad\'e method does not work can be easily constructed,
but we are unaware of a real problem for which
the solution is experimentally measured but the series can not be summarized
by the exciting combinations of Euler-McLaurin summation and
Pad\'e approximation.
That is why we believe that suggested method, criterion and software can be useful for many
applied problems in numerous areas not only in physics but in any scientific application where differential equations are solved.

\section{Acknowledgments} 
The authors are thankful to Evgeni Penev for the collaboration in early stages in the present research,
and for writing the $\ve$-algorithm in Fortran-90, Michail Mishonov gave a version in Java and 
Nedeltcho Zahariev and Zlatan Dimitrov in C-language.
The authors are also grateful to Claude Brezinski for pointing out and recommendation of 
Refs.~\cite{Gonnet:13,Beckermann:13}, to Ognyan Kounchev for
Refs.~\cite{Karlsson:76,Gonchar:78,Gonchar:89,Levis:18,Yattselev:18}, to Emil Horozov, Radostina Kamburova, Dantchi Koulova,
Peter Kenderov and Hassan Chamati for considerations, comments and stimulating discussions.
This work is partially supported by grant KP-06-N38/6 from 5.12.2019 of
the Bulgarian National Science Fund and National program
``Young scientists and postdoctoral researchers'' approved by DCM 577, 17.08.2018.

\appendix

\section{$\ve$-algorithm source code in Fortran-90}

This section includes the $\ve$-algorthm written in Fortran-90 subroutine and used for the calculations presented in the figures.
This Fortran program is our understaning of Cauhi-Jakobi problem
an if the time of floating point operation 
$\tau_{_\mathrm{FPO}}\rightarrow 0$
and also the machine epsilon 
$\epsilon_\mathrm{mach}\rightarrow 0$
we obtain the analytical result.
The program can be converted in a constructive proof.

\begin{lstlisting}[language=Fortran]
Subroutine LimesCross(N, S, rLimes, i_Pade, k_Pade, i_err, error, Nopt)
		implicit integer (i-n)
		implicit real(8) (a-h), real(8) (o-z)
		dimension s(0:N), a(0:N+2), b(0:N+1), c(0:N+1)
		Zero=0.d0
		One =1.d0
		rLimes=s(0)
		i_Pade=0
		k_Pade=0
                rNorm=Zero

		if (N<2) then
                	i_err=1 ! N>=2
                	Nopt=2
			rLimes=s(N)
			goto 137
		endif

		do i=0,N
			b(i) = s(i)
			c(i) = Zero
                        rNorm=rNorm+dabs(s(i))		
		end do
		if (rNorm==0) then 
                	i_err=2 ! S_i=0
			goto 137
		endif

                error= dabs(s(1)-s(0))
                do i=2, N
                eloc= dabs(s(i)-s(i-1))
                if(eloc<error) then
                	error = eloc
                	k_Pade=0
			i_Pade=i
			rLimes=s(i)
			i_err=3
                	Nopt=i_Pade + k_Pade+1 
                endif 
		end do
		
                Nb=N-2 ! N boundary 
                k=0 
123		continue

		!	beginning of the calculation of Pade table P_{k_Pade,i_Pade}
		!	k_Pade= power of the denominator
		!	i_Pade= power of the numerator
		!	Table ordered as a matrix numbering 
               
		do i=0, Nb
			Center= b(i+1)  
			if(k>=1) then
				rNorth= a(i+2)
				rNC= rNorth-Center
				if(rNC==Zero) then 
					rLimes= Center
					k_Pade=k
					i_Pade=k_Pade+i
					Nopt=i_Pade + k_Pade +1 
					error=zero
					i_err=4
					goto 137 
				endif
			endif
			West= b(i)
            WC= West-Center
            if(WC==Zero) then
               	rLimes= Center
               	k_Pade=k
               	i_Pade=k_Pade+i +1
               	Nopt=i_Pade + k_Pade +1 
               	error=zero
               	i_err=5
               	goto 137 
            endif
            East= b(i+2)
            EC= East-Center
            if(EC==Zero)  then 
               	rLimes= Center
               	k_Pade=k
               	i_Pade=k_Pade+i+2 
               	Nopt=i_Pade + k_Pade +1 
               	error=zero
               	i_err=6
               	goto 137 
            endif
!			The new criterion Wynn-Frobenius Eq. (11) from the present paper
			Wynn=One/WC + One/EC
			if(Wynn.NE.Zero) eloc=One/dabs(Wynn) ! in search of local minimum 
			denom=Wynn
			if(k>=1) denom=denom-One/rNC
			if(denom==Zero) then
				i_err=7 ! 1/rLimes=0, i.e. divergent result
				Nopt=i_Pade + k_Pade +1 
				goto 137 
			endif                          
			South= Center+One/denom
			C(i)=South
			SC=South-Center
!			Optimal Pade approximant                
            if(eloc<error) then
                error= eloc
				i_err=0
                k_Pade=k+1
                i_Pade=k_Pade+i 
                Nopt=i_Pade + k_Pade +1           	
                rLimes= South ! = P(k_Pade,i_Pade)
            endif      	
        end do ! i cycle  
        k=k+1
!       Change of the succesions 
        do i=0, Nb
			a(i)=b(i)
			b(i)=c(i)
		enddo
!		Reduction of the boundary           
        Nb=Nb-2
        if(Nb>=0) goto 123
!       Exit 
137	continue 
	end Subroutine LimesCross
\end{lstlisting}


\section{Aitken algorithm}
In order to explain the novelty of the present achievement 
in the beginning we will recall the well known Aitken algorithm.
Later on using this basis we will explain the general polynomial exponential.
We follow the well-known reference book by Bronstein and Semendyayev\cite{Bronshtein:55}.
If we need to calculate value of the function $\varphi_n(x)$ at some fixed argument $x$, 
one can use following schematics (``crest by crest'') which is especially convenient 
for computer programming
\begin{align}
x_0 - x \qquad &f_0 \nn \\
x_1 - x \qquad &f_1 \qquad (f_0,f_1) \nn  \\
x_2 - x \qquad &f_2 \qquad (f_0,f_2) \qquad (f_0,f_1,f_2) \nn  \\
\dots \dots \dots & \dots \dots \dots \nn  \\
x_n - x \qquad &f_n \qquad (f_0,f_n) \qquad (f_0,f_1,f_n) \dots (f_0,f_1, \dots ,f_n).\nn 
\end{align}
Every symbol $(f_0,\,f_1,\dots, f_k)$ denotes the value in the point $x$
of the interpolating polynomial function, build on nodes $x_0,\,x_1,\dots,x_k$,
i.e. $f_0=f(x_0)$, $f_1=f(x_1)$, \dots, $f_n=f(x_n).$
Those values calculates column by column, in the following manner.
Numbers of the column $(f_0,\,f_k)$ are obtained by
\be
(f_0,f_k) = \frac{(x_0-x) f_k - (x_k-x) f_0}{(x_0-x) - (x_k-x)}.
\ee
Every following column is calculating by the former one at the same schematics, for example
\be
(f_0,f_1,f_k) = \frac{(x_1-x) (f_0,f_k) - (x_k-x) (f_0,f_1)}{(x_1-x) - (x_k-x)}.
\label{second}
\ee
The rounding error is minimized if we reorder the points as
\be
|x_0-x|<|x_1-x|<|x_2-x|<\dots|x_n-x|.
\ee
In such a manner we obtain the series
\be
S_0\equiv f_0,\quad S_1\equiv(f_0,f_1),\quad S_2\equiv(f_0,f_1,f_2),
\quad \dots,\quad  S_n\equiv(f_0,f_1,f_2,\dots,f_n).
\ee
In order to obtain an approximation of limit of this sequence
\be
S=\lim_{n\rightarrow\infty}S_n
\ee
we use the Wynn algorithm CNEWS algorithm with minimal-$|\eta|$ criterion, 
which works even if the sequence $(S_0,\,S_1,\,S_2,\dots)$ is divergent.
It is only necessary the problem which we solve to be analytical one born 
from a real physical problem.

Let us repeat the Aitken formula as an algorithm.
Having fixed node point $x_0$ we can calculate the the linear approximation using a second
node point $\tilde x$
\be
\varphi(x;x_0,\tilde x)=f(\tilde x)+\frac{f_0-f(\tilde x)}{x_0-\tilde x}(x-\tilde x)
= \frac{(x_0-\tilde x) f(\tilde x) - (\tilde x-x) f_0}{(x_0-x) - (\tilde x-x)}
,
\ee
where node arguments $x_0$ and $\tilde x$ are parameters of the linear function.
In such a way omittin in the notations the argument $x$ we can calculate the sequence
\be
(f_0,f_1)=\varphi(x;x_0,x_1),\quad
(f_0,f_2)=\varphi(x;x_0,x_2),\quad
(f_0,f_3)=\varphi(x;x_0,x_3),\quad,\dots,\quad
(f_0,f_n)=\varphi(x;x_0,x_n).
\ee
In order to rewrite this conveniently for programming algorithm, we can
use the same notations for the changing new values
\be
f_1\leftarrow (f_0,f_1),\quad
f_2\leftarrow (f_0,f_2),\quad
f_3\leftarrow (f_0,f_3),\quad,\dots,\quad
f_k\leftarrow(f_0,f_n).
\ee
Then for the next column of the Aitken table according Eq.~(\ref{second}) 
using so re-denoted values we calculate 
\be
f_k\leftarrow \frac{(x_1-x) f_k - (x_k-x) f_1}{(x_1-x) - (x_k-x)}, 
\qquad\mbox{for}\quad k=2,\,3,\,\dots,n.
\label{second}
\ee
The formula for a general column $l$ is the same
\be
f_k\leftarrow \frac{(x_l-x) f_k - (x_k-x) f_l}{(x_l-x) - (x_k-x)}
=\frac{\left\vert\begin{array}{cc}
(x_l-x)&f_l\\
(x_k-x)&f_k
   \end{array}\right\vert}
{x_l-x_k}
, 
\qquad\mbox{for}\qquad k=l+1,\dots,n;
\label{second}
\ee
i.e. linearly interpolated value of the function at point $x$
is ascribed as the new function value at the node point $x_k$.
In other words polynomial interpolation is reduces to sequential linear
interpolations
This procedure repeats until $l=n-1$.
In such a way we program the series identically by the so defined polynomial interpolations 
\be
(f_0,\,f_1,\,f_2,\,\dots\,f_n)= (S_0,\,S_1,\,S_2,\,\dots\,S_n)
\ee 
using sequentially 
\be
x_0,\qquad (x_0,x_1),\qquad (x_0,x_1,x_2),\qquad (x_0,x_1,x_2,\dots,x_n)
\ee
set of nodes.
According to the general theory of Pad\'e approximants the Wynn CNEWS algorithm converts 
polynomial approximation into Pad\'e ones, and for $N=n-1$ we obtain $N$-point 
Pad\'e approximant.

If the first $K$ points coincide 
$x_0\rightarrow x_1,\dots,\rightarrow x_K$
the corresponding polynomial approximant is just the Taylor series
\be
\phi(x;x_0)=f_0+f_0^{\prime}(x-x_0)+\frac{1}{2}f_0^{\prime\prime}(x-x_0)^2
+\dots+\frac{1}{K!}f_0^{(K)}(x-x_0)^K,
\ee
where
$f_0^{\prime},\,f_0^{\prime\prime},\,\dots,\,f_0^{(K)}$ are the first $K$ derivatives
of the approximated function, and $x_0$ is the parameter of this series expansion.
Re-denoting
\be
f_0=\phi(x;x_0),\quad 
f_1=\phi(x;x_1),\quad 
f_2=\phi(x;x_2),\quad, \dots,\quad
f_n=\phi(x;x_n),
\ee
we can apply the regular Aitken procedure using the new Taylor expanded functional values.
In such a way finally we have solved the generalized 
Cauchy-Jacobi problem: to obtain a rational approximation of a function
passing through the N-points $x_1,\,x_2,\dots,x_N$
with fixed values of its first derivatives
$f_i^{\prime},\,f_i^{\prime\prime},\,\dots,\,f_i^{(K)}$ 
at every at those points $i=1,\dots,N.$
We gave only an idea for proof and a Fortran program illustrate our
understanding.

We will use this achievement for many applied problems.
For example, for solution of ordinary differential equation 
\be
\frac{\mathrm{d}y}{\mathrm{d}x}=\mathcal{F}(x,y).
\ee
Using former N-points with first derivative Taylor approximation $K=1$
we can calculate the functional value $y(x)$ using the explained method of Pad\'e extrapolation
and then calculate the first derivative at this point explicitly $f^\prime(x)=\mathcal{F}(x,y(x))$.
Renumbering later $x_{n+1}\leftarrow x_n$, $x_{n-1}\leftarrow x_n,$\dots
$x_2\leftarrow x_1,$ $x_1=x_0$ and $x_0\leftarrow x$ we can continue the procedure.
The number of former points is determined by the indices of the optimal Pad\'e approximants 
according to the minimal-$|\eta|$ criterion.
For the step $h$ for the next point $x=x_0+h$ one can apply many different criteria.
To our knowledge such combination of methods for solution of ordinary differential equation without fixed order is a new one.
We consider that content of this Appendix is implicitly given in the body of the article 
we give this clarification only for the beginner users.
\end{document}